\documentclass[11 pt, draftcls, onecolumn]{IEEEtran} % draftcls
\usepackage{amssymb,epsfig,color,cite,amsmath,amsfonts,mathrsfs,algorithm,algorithmic,tikz,bm}
\usepackage{epstopdf}
\usepackage{enumitem}

\newtheorem{lemma}{Lemma}[section]
\newtheorem{theorem}{Theorem}[section]

\newtheorem{remark}{Remark}[section]

\newtheorem{example}{Example}[section]

\newtheorem{assumption}{Assumption}[section]
\DeclareMathOperator{\sign}{Sgn}
\DeclareMathOperator{\diag}{diag}
\DeclareMathOperator{\argmax}{argmax}
\DeclareMathOperator{\argmin}{argmin}
\DeclareMathOperator{\rint}{rint}
\DeclareMathOperator{\gph}{gph}
\begin{document}
\title{Distributed Nonsmooth Optimization with Coupled Inequality Constraints via Modified Lagrangian Function}

%\author{Michael~Shell,~\IEEEmembership{Member,~IEEE,}
%        John~Doe,~\IEEEmembership{Fellow,~OSA,}
%        and~Jane~Doe,~\IEEEmembership{Life~Fellow,~IEEE}% <-this % stops a space
%\thanks{M. Shell was with the Department
%of Electrical and Computer Engineering, Georgia Institute of Technology, Atlanta,

\author{Shu~Liang,
        Xianlin~Zeng,
        and~Yiguang~Hong% <-this % stops a space
\thanks{S. Liang, X. Zeng and Y. Hong are with the Key Laboratory of Systems and Control,  Academy of Mathematics and Systems Science, Chinese Academy of Sciences, Beijing, 100190, China. e-mails: shuliang@amss.ac.cn (S. Liang), xianlin.zeng@amss.ac.cn (X. Zeng), yghong@iss.ac.cn (Y. Hong).}% <-this % stops a space
}

\maketitle

\begin{abstract}
This technical note considers a distributed convex optimization problem with nonsmooth cost functions and coupled nonlinear inequality constraints. To solve the problem, we first propose a modified Lagrangian function containing local multipliers and a
nonsmooth penalty function. Then we construct a distributed continuous-time algorithm by virtue of a projected primal-dual subgradient dynamics. Based on the nonsmooth analysis and Lyapunov function, we obtain the existence of the solution to the nonsmooth algorithm and its convergence.
\end{abstract}


% \fi
\IEEEpeerreviewmaketitle

\section{Introduction}
Distributed convex optimization has attracted intense research attention in recent years, due to its theoretic significance and broad applications in many research fields such as sensor networks, smart grids and social networks. 
Various models of distributed optimization have been proposed and studied in the literature. Most works have focused on consensus-based formulations, where each agent estimates the entire optimal solution via plentiful discrete-time algorithms (e.g., see \cite{Nedich2016Achieving, Zhu2012Distributed} and the references therein). Recently, more and more effort has also been done for distributed continuous-time algorithms (see \cite{Shi2013reaching,Gharesifard2014Distributed,Liu2015Second, Zeng2017Distributed} for instance), partly due to the development of its hardware implementation \cite{Forti2004Generalized} and flexible application in continuous-time physical systems \cite{Zhang2017Distributed}.

Here we consider distributed optimizations with separable cost functions and coupled constraints. In the presence of a coupled constraint, the feasible region of one agent's decision variable is influenced by some other agents' decision variables. If such a constraint is known by all the related agents, various algorithms were obtained, such as dual gradient algorithms \cite{Nedic2009Approximate, Necoara2015Linear}, primal-dual algorithms \cite{Nedic2009Subgradient, Feijer2010Stability,Cherukuri2016Asymptotic}, the saddle-point-like algorithm \cite{Niederlander2016Distributed}, and the distributed Newton-type algorithm \cite{Wei2013Distributed}.  However, coupled constraints may not be available to each agent in practice, and then the aforementioned algorithms may not work if there is no central coordinator in the network.  To deal with the challenges, \cite{Cherukuri2016Initialization,Yi2016Initialization} developed distributed initialization-free algorithms for the optimal resource allocation, while \cite{Zeng2016Continuous} proposed a distributed algorithm for the extended monotropic optimization. Note that \cite{Cherukuri2016Initialization,Yi2016Initialization,Zeng2016Continuous} considered coupled equality constraints.  Moreover, \cite{Chang2014Distributed} proposed a distributed algorithm for coupled inequality constraints, based on the average consensus technique to estimate the constraint functions along with a local primal-dual perturbed subgradient method. These distributed algorithms adopted local dynamics to evaluate the optimal dual solution instead of the original centralized one. On the other hand, all of them have to further employ auxiliary dynamics in order to guarantee the correctness and convergence, whereas the distributed design may become quite complicated in the case with coupled inequality constraints.

The objective of this note is to develop a distributed algorithm for nonsmooth convex optimization with coupled inequality constraints. We propose a modified Lagrangian function such that not only its saddle point yields the correct optimal solution to the original problem, but also its primal-dual subgradient dynamics is fully distributed. Particularly, we introduce local multipliers to decouple the constraints and employ a nonsmooth penalty function for the correctness. Based on this modified Lagrangian function, we propose a continuous-time projected subgradient algorithm for saddle-point computation. Our algorithm is fully distributed since each agent updates its local variables according to its local data and the information of its neighbors, without requiring any center in the network. Moreover, our algorithm only involves the primal variables and local multipliers, which yields a lower order dynamics than those in existing algorithms.

The rest organization is as follows: Section 2 provides necessary preliminaries, while Section 3 formulates the problem. Then Section 4 presents the main results to prove the convergence of our nonsmooth algorithm, and Section 5 gives two numerical examples. Finally, Section 6 gives some concluding remarks.

{\em Notations}: Denote $\mathbb{R}^n$ as the
$n$-dimensional real vector space and $\mathbb{R}^n_+$ as the positive orthant in $\mathbb{R}^n$. Denote $\bm{0}$ as a vector with each component being zero. For a vector $a\in \mathbb{R}^n$, $a\leq \bm{0}$ (or $a<\bm{0}$) means that each component of $a$ is less than or equal to zero (or less than zero). Denote $\|\cdot\|$ and $|\cdot|$ as the $\ell_2$-norm and $\ell_1$-norm for vectors, respectively. Denote $col(x_1,...,x_N) = (x_1^{T}, ... , x_n^{T})^{T}$ as the column vector stacked with column vectors $x_1,...,x_N$. For a set $\Omega\subset \mathbb{R}^n$, $\rint(\Omega)$ is the relative interior and
$d_\Omega(x) \triangleq \inf_{y\in \Omega}\|y-x\|$ is the distance function between point $x$ and set $\Omega$.
\section{Preliminaries}
In this section, we introduce relevant preliminary knowledge about convex analysis, differential inclusions, and graph theory.

A set $C \subseteq \mathbb{R}^n$ is {\em convex} if $\lambda z_1
+(1-\lambda)z_2\in C$ for any $z_1, z_2 \in C$ and $\lambda\in [0,\,1]$. For $x\in C$, the {\em tangent cone} to $C$ at $x$, denoted by $\mathcal{T}_C(x)$, is defined as
\begin{equation*}
\mathcal{T}_C(x) \triangleq \big\{\lim_{k\to\infty}\frac{x_k-x}{t_k}\,|\,x_k\in C, t_k>0,\text{ and }x_k \to x, t_k\to 0\big\},
\end{equation*}
while the {\em normal cone} to $C$ at $x$, denoted by $\mathcal{N}_C(x)$, is defined as
\begin{equation*}
\mathcal{N}_C(x) \triangleq \{v\in \mathbb{R}^n \,|\, v^T(y -x) \leq 0, \text{ for all } y\in C\}\text{.}
\end{equation*}
A projection operator is defined as $P_C(z) \triangleq \mathop{\argmin}_{x\in C}\|x-z\|$, and an operator that projects a point $z\in\mathbb{R}^n$ (or a set) onto the tangent cone $\mathcal{T}_C(x)$ is $\Pi_C(x;z) \triangleq P_{\mathcal{T}_C(x)}(z)$.

A function $f: C\to \mathbb{R}$ is said to be {\em convex} (or {\em strictly convex}) if $f(\lambda z_1
+(1-\lambda)z_2) \leq \text{ (or $<$) } \lambda f(z_1) + (1-\lambda)f(z_2)$ for any $z_1, z_2 \in C, z_1 \neq z_2$ and $\lambda\in (0,\,1)$. A function $g$ is said to be a {\em concave function} if $-g$ is a convex function.

A set-valued map $\mathcal{F}$ from $\mathbb{R}^n$ to $\mathbb{R}^n$ is a map associated with any $x\in \mathbb{R}^n$ a subset $\mathcal{F}(x)$ of $\mathbb{R}^n$. $\mathcal{F}$ is said to be {\em upper semicontinuous} at $x_0\in \mathbb{R}^n$ if for any open set $E$ containing $\mathcal{F}(x_0)$, there exists a neighborhood $D$ of $x_0$ such that $\mathcal{F}(D)\subset E$. We say that $\mathcal{F}$ is upper semicontinuous if it is so at every $x_0\in \mathbb{R}^n$. The graph of $\mathcal{F}$, denoted by $\gph\mathcal{F}$, is the set consists of all pairs $(x,y)$ satisfying $y\in \mathcal{F}(x)$.

A differential inclusion can be expressed as follows:
\begin{equation}\label{eq:DI}
\dot{x} \in \mathcal{F}(x), \quad x(0) = x_0.
\end{equation}
A map $x(t): [0, +\infty) \rightarrow \mathbb{R}^n$ is said to be a {\em solution} to \eqref{eq:DI} if it is  absolutely continuous and satisfies the inclusion for almost all $t \in [0, +\infty)$.

A graph of a network is denoted by $\mathcal{G}=(\mathcal{V},\,\mathcal{E})$, where $\mathcal{V} = \{1,...,N\}$ is a set of nodes and $\mathcal{E} \subseteq
\mathcal{V}\times \mathcal{V} $ is a set of edges. Node $j$ is said to be a {\em neighbor} of node $i$ if $\{i,j\}\in \mathcal{E}$. The set of all the neighbors of node $i$ is denoted by $\mathcal{N}_i$. $\mathcal{G}$ is said to be {\em undirected} if $(i,\,j)\in \mathcal{E} \Leftrightarrow (j,\,i)\in \mathcal{E}$. A path of $\mathcal{G}$ is a sequence of distinct nodes where any pair of consecutive nodes in the sequence
has an edge of $\mathcal{G}$. Node $j$ is said
to be {\em connected} to node $i$ if there is a path from $j$ to $i$. $\mathcal{G}$ is said to be connected if any two nodes are
connected.  The detailed knowledge about graph theory can be found in \cite{Godsil01}.

The following lemma collects some results given in \cite{Clarke1998Nonsmooth} that will be used in our analysis.%(Proposition 1.5, page 73, Corollary 2.7 page 37, and Proposition 6.3, page 50), {Theorem 3.33, page 126}.
\begin{lemma}\label{lem:subdifferential}
Let $f: X \to \mathbb{R}$ be locally Lipschitz continuous and let $\Omega \subset X$ be a closed convex subset, where $X\subset \mathbb{R}^n$. Then the following statements hold.
\begin{enumerate}[leftmargin=*]
\item \label{item:subdifferential2} $\gph\partial f$ is closed.
\item \label{item:subdifferential3} If $f$ is convex, then $x^*\in \mathop{\argmin}_{x\in \Omega}f(x)$ if and only if $0\in \partial f(x^*) + \mathcal{N}_\Omega(x^*)$.
\item \label{item:subdifferential4} Suppose $f$ has a Lipschitz constant $K_0$ on an open set that contains $\Omega$. When $K>K_0$, $x^*\in \mathop{\argmin}_{x\in \Omega}f(x)$ if and only if $x^*\in \mathop{\argmin}_{x\in X}f(x)+Kd_\Omega(x)$.
\end{enumerate}
\end{lemma}

A collection of results in \cite{Aubin1984Differential} with respect to set-valued maps and differential inclusions are given below.
\begin{lemma}\label{lem:setvalue}
The following statements hold.
\begin{enumerate}[leftmargin=*]
\item \label{item:setvalue1}A set-valued map $\mathcal{F}$ from $\mathbb{R}^n$ to $\mathbb{R}^n$ is upper semicontinuous if it has compact values and $\gph\mathcal{F}$ is closed.

\item \label{item:setvalue3}Let $\mathcal{F}_0$ from $\mathbb{R}^n$ to $\mathbb{R}^n$ be a set-valued map and $C\subset \mathbb{R}^n$ be a closed convex subset. Consider the following two differential inclusions
\begin{align}
\label{eq:di1}
\dot x &\in \mathcal{F}_0(x) - \mathcal{N}_{C}(x), \quad x(0) = x_0\in C,\\
\label{eq:di2}
\dot x & \in \Pi_{C}(x,\mathcal{F}_0(x)), \quad x(0) = x_0 \in C.
\end{align}
Then $x(\cdot)$ is a solution to \eqref{eq:di1} if and only if it is a solution to \eqref{eq:di2}.

\item \label{item:setvalue4} For any $x_0\in C$, there is a solution to the differential inclusion \eqref{eq:di1} if $\mathcal{F}_0$ is upper semicontinuous and $C$ is compact and convex.
\end{enumerate}
\end{lemma}

Moreover, we introduce a lemma from \cite[Lemma 4.2]{Zeng2017Distributed} which will be used in the convergence analysis.

\begin{lemma}\label{lem:semistability}
Let $x(t)$ be a solution to the differential inclusion \eqref{eq:DI}. If $z$ is a Lyapunov stable equilibrium of \eqref{eq:DI} and is also a cluster point of $x(\cdot)$, then $\lim_{t\to+\infty}x(t) = z$.
\end{lemma}
\section{Problem Formulation}

Consider a multi-agent network with $N$ agents, whose label set is denoted as $\mathcal{V} = \{1,...,N\}$, cooperating over a graph $\mathcal{G} = \{\mathcal{V},\mathcal{E}\}$. For each agent $i$, there are a local decision variable $x_i\in \mathbb{R}^{n_i}$ and a local constraint set $\Omega_i\subset \mathbb{R}^{n_i}$ for $i\in \mathcal{V}$. Define $\bm{x}\triangleq col(x_1,...,x_N)$, and define
the total cost function of the network as $f(\bm{x}) \triangleq \sum_{i\in \mathcal{V}}f_i(x_i)$, where $f_i:\Omega_i\to \mathbb{R}$ is a (nonsmooth) local cost function of agent $i$. In addition, the agents are subject to coupled inequality constraints in the form of $\bm{g}(\bm{x}) \triangleq \sum_{i\in\mathcal{V}} \bm{g}_i(x_i)\leq \bm{0}$, where $\bm{g}_i:\Omega_i\to\mathbb{R}^M$ are continuous mappings for $i\in \mathcal{V}$ (that is, $\bm{g}_i = (g_{i1},...,g_{iM})^T$ and $g_{ik}:\Omega_i \to \mathbb{R}$ are continuous functions for all $i$ and $k$). To be strict, the optimization problem can be formulated as:
\begin{equation}\label{eq:optimizationProblem}
\min_{\bm{x}\in\Omega} f(\bm{x}), \text{ s.t. } \bm{g}(\bm{x})\leq \bm{0},
\end{equation}
where $\Omega\triangleq \prod_{i\in \mathcal{V}}\Omega_i \subset \mathbb{R}^{n_1+\cdots+n_N}$ denotes the local constraints of $N$ agents.

The following assumption is needed to ensure the well-posedness of problem \eqref{eq:optimizationProblem}.
\begin{assumption}\label{assum:1}
~
\begin{enumerate}[leftmargin=*]
\item (Convexity and continuity) For all $i\in\mathcal{V}$, $\Omega_i$ is compact and convex. On an open set containing $\Omega_i$, $f_i$ is strictly convex and $\bm{g}_i$ is convex, and $f_i$ and $\bm{g}_i$ are locally Lipschitz continuous.
\item (Slater's constraint qualification) There exists $\bar{\bm{x}}\in \rint(\Omega)$ such that $\bm{g}(\bar{\bm{x}})< \bm{0}$.
\item (Communication topology) The graph $\mathcal{G}$ is connected and undirected.
\end{enumerate}
\end{assumption}

This assumption is quite mild and similar ones are widely used in the literature (e.g., \cite{Yi2016Initialization}).

The goal of this note is to develop a {\em distributed continuous-time algorithm} for solving \eqref{eq:optimizationProblem} with each agent communicating with their neighbors. Moreover, for every $i\in\mathcal{V}$, agent $i$ can only access $\bm{g}_i(x_i)$ rather than $\bm{g}(\bm{x})$.

The differences between our problem and those in existing literature are as follows.
\begin{itemize}[leftmargin=*]
\item The decision variables can be heterogeneous with possibly different dimensions, in contrast to those consensus-based models.
\item The cost and constraint functions can be nonsmooth, while some projected dynamics \cite{Feijer2010Stability,Cherukuri2016Asymptotic,Yi2016Initialization} and Newton type method \cite{Wei2013Distributed} depend on smoothness.
\item The coupled constraints may be unavailable to local agents, different from \cite{Niederlander2016Distributed}. Also, coupled inequality constraints are considered as in \cite{Chang2014Distributed}, different from the coupled (affine) equality ones studied in \cite{Cherukuri2016Initialization,Yi2016Initialization,Zeng2016Continuous}.
\end{itemize}

\section{Main results}
In this section, we first propose a modified Lagrangian function and then propose a distributed continuous-time algorithm for the considered optimization problem. Moreover, we prove the existence of the solution to the nonsmooth algorithm along with the discussions on its convergence and the rate.

\subsection{Lagrangian Function and Distributed Algorithm Design}
Consider the following dual problem with respect to the primal one \eqref{eq:optimizationProblem},
\begin{equation}\label{eq:dualProblem}
\max_{\lambda\geq \bm{0}}q(\lambda), \quad q(\lambda) \triangleq \min_{\bm{x}\in \Omega}\mathcal{L}(\bm{x},\lambda),
\end{equation}
where $\mathcal{L}:\Omega\times \mathbb{R}^M_+ \to \mathbb{R}$ is the Lagrangian function defined as
\begin{equation}\label{eq:Lagrangian}
\mathcal{L}(\bm{x},\lambda) \triangleq f(\bm{x}) + \lambda^T \bm{g}(\bm{x}).
\end{equation}
It has been shown in \cite{Nedic2009Approximate} that the optimal dual solution $\lambda^*$ of \eqref{eq:dualProblem} lies in a compact set $\mathcal{D}\subset \mathbb{R}^M_+$, given by
\begin{equation}\label{eq:dualSet}
\mathcal{D} \triangleq \{\lambda\in \mathbb{R}^M_+\,|\, \|\lambda\|\leq \frac{f(\bar{\bm{x}}) - \tilde{q}}{\gamma}\},
\end{equation}
where $\bar{\bm{x}}$ is a Slater point of \eqref{eq:optimizationProblem}, $\tilde{q} = \min_{\bm{x}\in \Omega}\mathcal{L}(\bm{x},\tilde{\lambda})$ is a dual function value for an arbitrary $\tilde{\lambda}\geq 0$, $\gamma = \min_{k = 1,...,M}\{-\sum_{i\in\mathcal{V}}g_{ik}(\bar{x}_i)\}$.

We present the following lemma, which is a well-known convex optimization result \cite{Rockafellar1998Variational}.

\begin{lemma}\label{lem:collection}
Under Assumption \ref{assum:1}, the following statements are equivalent:
\begin{enumerate}[leftmargin=*]
\item (Primal-dual characterization) $(\bm{x}^*, \lambda^*)$ is a  primal-dual solution pair of problems \eqref{eq:optimizationProblem} and \eqref{eq:dualProblem}.
\item (Saddle-point characterization) $(\bm{x}^*, \lambda^*)$ is a {\em saddle point} of Lagrangian function \eqref{eq:Lagrangian}, that is,
\begin{equation*}
\mathcal{L}(\bm{x}^*,\lambda) \leq \mathcal{L}(\bm{x}^*,\lambda^*) \leq \mathcal{L}(\bm{x},\lambda^*), \quad \forall\, \bm{x}\in \Omega,\,\lambda\geq \bm{0}.
\end{equation*}
\item (KKT characterization) $(\bm{x}^*, \lambda^*)$ satisfies
\begin{equation*}
\bm{0} \in \partial f(\bm{x}^*) + \partial\bm{g}(\bm{x}^*)\lambda^{*} + \mathcal{N}_{\Omega}(\bm{x}^*), \quad \bm{0} \leq \lambda^* \perp -\bm{g}(\bm{x}^*) \geq \bm{0}.
\end{equation*}
\item (Minimax characterization) $(\bm{x}^*,\lambda^*)$ is a solution of the minimax problem
\begin{equation*}
\min_{x\in \Omega}\{\max_{\lambda\in \mathcal{D}} \mathcal{L}(\bm{x},\lambda)\}.
\end{equation*}
\end{enumerate}
Moreover, since $f(\bm{x})$ is strictly convex, such $\bm{x}^*$ is unique while $\lambda^*$ may not be unique.
\end{lemma}

A centralized projected primal-dual algorithm with respect to $\mathcal{L}(\bm{x},\lambda)$ can be written as (referring to \cite{Cherukuri2016Asymptotic}):
\begin{equation*}
\left\{\begin{aligned}
\dot{\bm{x}} &\in \Pi_{\Omega}(\bm{x}, - \partial f(\bm{x}) - \partial\bm{g}(\bm{x})\lambda), \, &&\bm{x}(0) \in \Omega,\\
\dot{\lambda} &\in \Pi_{\mathbb{R}^{M}_+}(\lambda, \bm{g}_1(x_1) + \cdots + \bm{g}_N(x_N)), \, &&\lambda(0) \in \mathbb{R}^{M}_+,
\end{aligned}\right.
\end{equation*}
which needs a center to broadcast $\lambda$ and gather $\bm{g}_1, ..., \bm{g}_N$ for the update. In order to develop fully distributed algorithms without a center, we employ local multipliers and a nonsmooth penalty function to construct a modified Lagrangian function. To be specific, define
\begin{subequations}
\begin{align}
\bm{\lambda} & \triangleq col(\lambda_1,..., \lambda_N)\in \mathbb{R}^{MN}_+,\\
\label{eq:S}
\mathcal{S} & \triangleq \{\bm{\lambda}\in \mathbb{R}^{MN}_+\,|\, \lambda_1 = \cdots = \lambda_N\},\\
\label{eq:phi}
\phi(\bm{\lambda}) & \triangleq \frac{1}{2} \sum_{i\in \mathcal{V}} \sum_{j\in \mathcal{N}_i}|\lambda_i-\lambda_j|\text{,}\\
\label{eq:tildeL}
\tilde{\mathcal{L}}(\bm{x},\bm{\lambda}) & \triangleq \sum_{i\in \mathcal{V}}f_i(x_i)+\lambda_i^T\bm{g}_i(x_i) - K\phi(\bm{\lambda}),
\end{align}
\end{subequations}
where $\bm{\lambda}$ is a collection of {\em local multipliers} employed for distributed design and $\mathcal{S}$ is a cone for the local multipliers to reach a consensus there.  Moreover, $\phi(\bm{\lambda})$ serves as a {\em metric of consensus} for the multipliers and $\tilde{\mathcal{L}}(\bm{x},\bm{\lambda})$ is a modified Lagrangian function with a constant $K>0$.

The following lemma reveals that the nonsmooth $K\phi(\bm{\lambda})$ plays a role as an exact penalty function for the consensus of multipliers.

\begin{lemma}\label{lem:exactPenalty}
Under Assumption \ref{assum:1}, for any $\bm{x}\in \Omega$, there holds
\begin{equation*}
\mathop{\argmax}_{\bm{\lambda}\in \mathbb{R}^{MN}_+} \tilde{\mathcal{L}}(\bm{x},\bm{\lambda}) = \mathop{\argmax}_{\bm{\lambda}\in \mathcal{S}} \tilde{\mathcal{L}}(\bm{x},\bm{\lambda}),
\end{equation*}
provided $K > \sqrt{N}K_0$, where $K_0\triangleq \max_{\bm{x}\in \Omega}\|col(\bm{g}_1(x_1),...,\bm{g}_N(x_N))\|$.
\end{lemma}

\begin{IEEEproof}
It follows from part \ref{item:subdifferential4}) in Lemma \ref{lem:subdifferential} that, for any $K_d>K_0$,
\begin{equation*}
\mathop{\argmax}_{\bm{\lambda}\in \mathcal{S}}\sum_{i\in \mathcal{V}}f_i(x_i)+\lambda_i^T\bm{g}_i(x_i) = \mathop{\argmax}_{\bm{\lambda}\in \mathbb{R}^{MN}_+} \sum_{i\in \mathcal{V}}f_i(x_i)+\lambda_i^T\bm{g}_i(x_i) - K_dd_{\mathcal{S}}(\bm{\lambda}).
\end{equation*}
Since $\phi(\bm{\lambda}) = d_{\mathcal{S}}(\bm{\lambda}) = 0, \forall\, \bm{\lambda} \in \mathcal{S}$, it suffices to prove $\sqrt{N}\phi(\bm{\lambda}) > d_{\mathcal{S}}(\bm{\lambda})$ for all $\bm{\lambda}\in \mathbb{R}^{MN}_+\setminus \mathcal{S}$.

On one hand,
\begin{equation*}
d_{\mathcal{S}}^2(\bm{\lambda}) = \min_{\tilde{\bm{\lambda}}\in S} \|\tilde{\bm{\lambda}} - \bm{\lambda}\|^2= \sum_{k=1}^N \big\|\lambda_k - \frac{\lambda_1 + \cdots + \lambda_N}{N}\big\|^2 \leq \frac{1}{N}\sum_{k=1}^N \sum_{l=1}^N\|\lambda_k-\lambda_l\|^2 \leq \frac{1}{N}\sum_{k=1}^N \sum_{l=1}^N|\lambda_k-\lambda_l|^2.
\end{equation*}
On the other hand, since the graph is connected and undirected, there is a path $\mathcal{P}_{kl}\subset \mathcal{E}$ connecting nodes $k$ and $l$ for any $k,l\in\mathcal{V}$. Then
\begin{align*}
\phi(\bm{\lambda}) = \frac{1}{2}\sum_{(i,j)\in \mathcal{E}}|\lambda_i - \lambda_j| \geq \frac{1}{2}\sum_{(i,j)\in \mathcal{P}_{kl}}|\lambda_i - \lambda_j| \geq |\lambda_k-\lambda_l|.
\end{align*}
Thus, $d_{\mathcal{S}}^2(\bm{\lambda}) \leq N\phi^2(\bm{\lambda})$ and the equality holds if and only if $\bm{\lambda} \in \mathcal{S}$, which implies the conclusion.
\end{IEEEproof}

The correctness of the Lagrangian function $\tilde{\mathcal{L}}(\bm{x},\bm{\lambda})$ to problem \eqref{eq:optimizationProblem}  is indicated in the following result.

\begin{theorem}\label{thm:equilibrium}
Under Assumption \ref{assum:1}, the following statements are equivalent:
\begin{enumerate}[leftmargin=*]
\item $(\bm{x}^*, \bm{\lambda}^*) \in \Omega \times \mathbb{R}^{MN}_+$ renders the following equations
\begin{subequations}\label{eq:equilibrium}
\begin{align}
0 &\in \Pi_{\Omega}(\bm{x}^*, - \partial_{\bm{x}}\tilde{\mathcal{L}}(\bm{x}^*,\bm{\lambda}^*)),\\
0 &\in \Pi_{\mathbb{R}^{MN}_+}(\bm{\lambda}^*, -\partial_{\bm{\lambda}}(-\tilde{\mathcal{L}})(\bm{x}^*,\bm{\lambda}^*)).
\end{align}
\end{subequations}

\item $(\bm{x}^*, \bm{\lambda}^*)$ is a saddle point of $\tilde{\mathcal{L}}(\bm{x},\bm{\lambda})$ in $\Omega \times \mathbb{R}^{MN}_+$.

\item $\bm{\lambda^*} = col(\lambda^*,...,\lambda^*)$ and $(\bm{x}^*, \lambda^*)$ is a saddle point of $\mathcal{L}(\bm{x},\lambda)$ in $\Omega \times \mathbb{R}^{M}_+$.
\end{enumerate}
\end{theorem}

\begin{IEEEproof}
1) $\Rightarrow$ 2): Let $(\bm{x}^*, \bm{\lambda^*})\in \Omega \times \mathbb{R}^{MN}_+$ satisfying \eqref{eq:equilibrium}. Then
\begin{subequations}\label{eq:equlibriumCondition}
\begin{align}
\label{eq:equlibriumCondition1}
0 &\in - \partial_{\bm{x}}\tilde{\mathcal{L}}(\bm{x}^*,\bm{\lambda}^*) - \mathcal{N}_\Omega(\bm{x}^*),\\
\label{eq:equlibriumCondition2}
0 &\in -\partial_{\bm{\lambda}}(-\tilde{\mathcal{L}})(\bm{x}^*,\bm{\lambda}^*) - \mathcal{N}_{\mathbb{R}^{MN}_+}(\bm{\lambda}^*).
\end{align}
\end{subequations}
Since $\tilde{\mathcal{L}}(\bm{x},\bm{\lambda})$ is convex in $\bm{x}$ and concave in $\bm{\lambda}$ (or equivalently, $-\tilde{\mathcal{L}}(\bm{x},\bm{\lambda})$ is convex in $\bm{\lambda}$), it follows from part \ref{item:subdifferential3}) in Lemma \ref{lem:subdifferential} that $\bm{x}^*$ is the minimum point of $\tilde{\mathcal{L}}(\cdot,\bm{\lambda}^*)$ in $\Omega$ and $\bm{\lambda}^*$ is a maximum point of $\tilde{\mathcal{L}}(\bm{x}^*,\cdot)$ in $\mathbb{R}^{MN}_+$, which implies statement 2).

2) $\Rightarrow$ 3): Let $(\bm{x}^*, \bm{\lambda}^*)$ be a saddle point of $\tilde{\mathcal{L}}(\bm{x},\bm{\lambda})$ in $\Omega \times \mathbb{R}^{MN}_+$. Then
\begin{equation}\label{eq:equlibriumCondition3}
\bm{\lambda}^* \in \mathop{\argmax}_{\bm{\lambda}\geq \bm{0}} \tilde{\mathcal{L}}(\bm{x}^*,\bm{\lambda}).
\end{equation}
It follows from Lemma \ref{lem:exactPenalty} that $\bm{\lambda}^* = col(\lambda^*,...,\lambda^*)$. Substituting this $\bm{\lambda}^*$ into the saddle point inequalities with respect to $\tilde{\mathcal{L}}(\bm{x},\bm{\lambda})$ yields
\begin{equation}\label{eq:equlibriumCondition4}
\bm{x}^* \in \mathop{\argmin}_{\bm{x}\in \Omega}\mathcal{L}(\bm{x},\lambda^*)\text{ and } \lambda^* \in \mathop{\argmax}_{\lambda\geq \bm{0}} \mathcal{L}(\bm{x}^*,\lambda),
\end{equation}
because of the identity $\tilde{\mathcal{L}}(\bm{x},\bm{\lambda}^*) = \mathcal{L}(\bm{x},\lambda^*)$. Therefore, the conclusion follows.

3) $\Rightarrow$ 1): Suppose $\bm{\lambda^*} = col(\lambda^*,...,\lambda^*)$ and $(\bm{x}^*, \lambda^*)$ is a saddle point of $\mathcal{L}(\bm{x},\lambda)$. According to Lemma \ref{lem:exactPenalty}, condition \eqref{eq:equlibriumCondition3} holds. Again, from part \ref{item:subdifferential3}) in Lemma \ref{lem:subdifferential}, condition \eqref{eq:equlibriumCondition} holds, which implies statement 1).
\end{IEEEproof}

By Theorem \ref{thm:equilibrium} and Lemma \ref{lem:collection}, the saddle points of $\tilde{\mathcal{L}}(\bm{x},\bm{\lambda})$ match exactly the saddle points of $\mathcal{L}(\bm{x},\lambda)$, which are in accordance with the optimal primal-dual solutions.

Based on $\tilde{\mathcal{L}}(\bm{x},\bm{\lambda})$, we present a distributed continuous-time algorithm to solve \eqref{eq:optimizationProblem} as follows:
\begin{equation}\label{eq:distributedAlgorithm}
\forall\,i\in\mathcal{V}:\, \left\{\begin{aligned}
\dot{x}_i &\in \Pi_{\Omega_i}(x_i, - \partial_{x_i}\tilde{\mathcal{L}}(\bm{x},\bm{\lambda})), \, && x_i(0) \in \Omega_i\\
\dot{\lambda}_i &\in \Pi_{\mathbb{R}^M_+}(\lambda_i, -\partial_{\lambda_i}(-\tilde{\mathcal{L}})(\bm{x},\bm{\lambda})), \, && \lambda_i(0) \in \mathbb{R}^M_+,
\end{aligned}\right.
\end{equation}
where
\begin{equation*}
\partial_{x_i}\tilde{\mathcal{L}}(\bm{x},\bm{\lambda}) =  \partial f_i(x_i) + \partial \bm{g}_{i}(x_i)\lambda_i, \quad -\partial_{\lambda_i}(-\tilde{\mathcal{L}})(\bm{x},\bm{\lambda}) = \bm{g}_i(x_i) - K\sum_{j\in \mathcal{N}_i}\sign(\lambda_i-\lambda_j),
\end{equation*}
(the second equality follows because graph $\mathcal{G}$ is undirected) and $\sign(\cdot)$ is the set-valued sign function with each component defined as
\begin{equation*}
\sign(y) \triangleq \partial |y| = \left\{\begin{aligned}
& \{1\} &&\text{  if } y>0\\
& \{-1\} &&\text{  if } y<0\\
& [-1,\,1] &&\text{  if } y=0\\
\end{aligned}\right.\text{.}
\end{equation*}
For simplicity, we rewrite algorithm \eqref{eq:distributedAlgorithm} in a compact form
\begin{equation}\label{eq:distributedAlgorithmCompact}
\left\{\begin{aligned}
\dot{\bm{x}} &\in \Pi_{\Omega}(\bm{x}, - \partial_{\bm{x}}\tilde{\mathcal{L}}(\bm{x},\bm{\lambda})), \, && \bm{x}(0) \in \Omega\\
\dot{\bm{\lambda}} &\in \Pi_{\mathbb{R}^{MN}_+}(\bm{\lambda}, -\partial_{\bm{\lambda}}(-\tilde{\mathcal{L}})(\bm{x},\bm{\lambda})), \, && \bm{\lambda}(0) \in \mathbb{R}^{MN}_+.
\end{aligned}\right.
\end{equation}

Algorithm \eqref{eq:distributedAlgorithm} is fully distributed since each agent $i\in\mathcal{V}$ only updates its local variables $x_i$ and $\lambda_i$ according to its local functions $f_i, \bm{g}_i$ and the information of its neighbors $\lambda_j, j\in \mathcal{N}_i$.
\begin{remark}
Some discussions about our method are given below.
\begin{itemize}[leftmargin=*]
\item In the original $\mathcal{L}(\bm{x},\lambda)$, each $\bm{g}_i(x_i)$ shares a common $\lambda$ and the multiplier $\lambda$ performs on the coupled $\bm{g}(\bm{x})$, while, in the modified one $\tilde{\mathcal{L}}(\bm{x},\bm{\lambda})$, all the local parts of cost and constraint functions are gathered in a decoupled way.

\item Our distributed algorithm involves only primal variables and local multipliers without auxiliary dynamics, while some existing distributed algorithms such as those given in \cite{Cherukuri2016Initialization,Yi2016Initialization,Zeng2016Continuous,Chang2014Distributed} employed auxiliary dynamics for the convergence.

\item From the estimation $K_0 \leq \sum_{i=1}^N \max_{x_i \in \Omega_i}\{\|g_i(x_i)\|\}$ for Lemma \ref{lem:exactPenalty}, parameter $K$ can be assigned via local estimation of each $\max_{x_i \in \Omega_i}\|g_i(x_i)\|$ and calculating the sum in a distributed manner.
\end{itemize}
\end{remark}

\subsection{Existence and Convergence}
Dynamics \eqref{eq:distributedAlgorithmCompact} is nonsmooth due to the projection operator and subgradients of the nonsmooth Lagrangian function. Thus, we need to check the existence of its solution (trajectory).

\begin{theorem}\label{thm:existence}
Under Assumption \ref{assum:1}, for any initial value $\bm{x}(0)\in \Omega, \bm{\lambda}(0) \in \mathbb{R}^{MN}_+$, there exists a solution to \eqref{eq:distributedAlgorithmCompact}.
\end{theorem}

\begin{IEEEproof}
Let $\tilde{\mathcal{D}}$ be the convex hull of $\bm{\lambda}(0)$ and $\prod_{i=1}^N\mathcal{D}$, where $\mathcal{D}$ is in \eqref{eq:dualSet}. Then $\tilde{\mathcal{D}}$ is compact and convex. Consider the following differential inclusion
\begin{equation}\label{eq:compactDI}
\left\{\begin{aligned}
\dot{\bm{x}} &\in \Pi_{\Omega}(\bm{x}, - \partial_{\bm{x}}\tilde{\mathcal{L}}(\bm{x},\bm{\lambda}))\\
\dot{\bm{\lambda}} &\in \Pi_{\tilde{\mathcal{D}}}(\bm{\lambda}, -\partial_{\bm{\lambda}}(-\tilde{\mathcal{L}})(\bm{x},\bm{\lambda}))
\end{aligned}\right.\text{.}
\end{equation}
Since $\mathcal{T}_{\tilde{\mathcal{D}}}(\bm{\lambda}) \subset \mathcal{T}_{\mathbb{R}^{MN}_+}(\bm{\lambda}), \forall\, \bm{\lambda}\in \tilde{\mathcal{D}}$, any solution to \eqref{eq:compactDI} is also a solution to \eqref{eq:distributedAlgorithmCompact}. Thus, it suffices to prove the existence of solution for \eqref{eq:compactDI}.

Let $\mathcal{F}(\bm{x},\bm{\lambda}) \triangleq col(- \partial_{\bm{x}}\tilde{\mathcal{L}}(\bm{x},\bm{\lambda}), -\partial_{\bm{\lambda}}(-\tilde{\mathcal{L}})(\bm{x},\bm{\lambda}))$, and $C \triangleq \Omega\times \tilde{\mathcal{D}}$. We claim that $\mathcal{F}$ is upper semicontinuous over $C$. The locally Lipschtiz continuity of $f(\bm{x}), \bm{g}(\bm{x}), \phi(\bm{\lambda})$ implies that $\mathcal{F}$ has compact values over the compact set $C$. Then it suffices to prove $\gph\mathcal{F}$ is closed due to part \ref{item:setvalue1}) in Lemma \ref{lem:setvalue}. Let $\{\bm{x}_k, \bm{\lambda}_k\}$ and $\{\zeta_k,\eta_k\}$ be sequences in $C$ and $\mathbb{R}^{n_1+\cdots+n_N}\times \mathbb{R}^{MN}$ such that (a) $col(\zeta_k,\eta_k) \in \mathcal{F}(\bm{x}_k,\bm{\lambda}_k)$, (b) $(\bm{x}_k, \bm{\lambda}_k)$ converges to $(\bm{x}, \bm{\lambda})$, and (c) $(\zeta,\eta)$ is a cluster point of the sequence $(\zeta_k,\eta_k)$. We can extract a subsequence of $(\zeta_k,\eta_k)$ (without relabeling) such that $\lim_{k\to+\infty}(\zeta_k,\eta_k) = (\zeta,\eta)$. Since $- \partial_{\bm{x}}\tilde{\mathcal{L}}(\bm{x},\bm{\lambda}) = - \partial f(\bm{x}) - \diag\{\partial \bm{g}_1(\bm{x}), ..., \partial \bm{g}_N(\bm{x})\}\bm{\lambda}$,  $\zeta_k = -\alpha_k - \beta_k\bm{\lambda}_k$, where $\alpha_k\in \partial f(\bm{x}_k)$ and $\beta_k \in \diag\{\partial \bm{g}_1(\bm{x}_k), ..., \partial \bm{g}_N(\bm{x}_k)\}$. It follows from part \ref{item:subdifferential2}) in Lemma \ref{lem:subdifferential} that $\lim_{k\to+\infty} \alpha_k = \alpha \in \partial f(\bm{x})$ and $\lim_{k\to+\infty} \beta_k = \beta \in \diag\{\partial \bm{g}_1(\bm{x}), ..., \partial \bm{g}_N(\bm{x})\}$ after extracting subsequences of $\{\alpha_k\}$ and $\{\beta_k\}$ without relabeling. Therefore, $\zeta \in - \partial_{\bm{x}}\tilde{\mathcal{L}}(\bm{x},\bm{\lambda})$. Similarly, $\eta \in -\partial_{\bm{\lambda}}(-\tilde{\mathcal{L}})(\bm{x},\bm{\lambda})$. Thus, $(\zeta,\eta) \in \mathcal{F}(\bm{x},\bm{\lambda})$, i.e., $\gph\mathcal{F}$ is closed.

Finally, according to part \ref{item:setvalue3}) and part \ref{item:setvalue4}) in Lemma \ref{lem:setvalue}, there exists a solution to system \eqref{eq:compactDI}, which is also a solution to \eqref{eq:distributedAlgorithmCompact}.
\end{IEEEproof}

Then it is time to show the convergence of our algorithm.

\begin{theorem}\label{thm:convergence}
Under Assumption \ref{assum:1}, algorithm \eqref{eq:distributedAlgorithmCompact} is stable and any of its solutions converges to the set of saddle points of $\tilde{\mathcal{L}}$. Moreover, for any solution $(\bm{x}(t),\bm{\lambda}(t))$, there exists a saddle point $(\bm{x}^*,\tilde{\bm{\lambda}}^*)$ of $\tilde{\mathcal{L}}$ such that
\begin{equation}\label{eq:convergence}
\lim_{t\to+\infty} (\bm{x}(t),\bm{\lambda}(t)) = (\bm{x}^*,\tilde{\bm{\lambda}}^*).
\end{equation}
\end{theorem}

\begin{IEEEproof}
For all $\bm{x}, \underline{\bm{x}}\in \Omega, \bm{\lambda}, \underline{\bm{\lambda}} \in \mathbb{R}_+^{MN}$, the following basic conditions hold, according to the definitions of projection, normal cone and the convexity-concavity of $\tilde{\mathcal{L}}$.
\begin{itemize}[leftmargin=*]
\item (projection)
\begin{align*}
&\Pi_{\Omega}(\bm{x}, - \partial_{\bm{x}}\tilde{\mathcal{L}}(\bm{x},\bm{\lambda})) \subset - \partial_{\bm{x}}\tilde{\mathcal{L}}(\bm{x},\bm{\lambda}) - \mathcal{N}_\Omega(\bm{x}), \\ &\Pi_{\mathbb{R}^{MN}_+}(\bm{\lambda}, -\partial_{\bm{\lambda}}(-\tilde{\mathcal{L}})(\bm{x},\bm{\lambda})) \subset -\partial_{\bm{\lambda}}(-\tilde{\mathcal{L}})(\bm{x},\bm{\lambda}) - \mathcal{N}_{\mathbb{R}^{MN}_+}(\bm{\lambda}),
\end{align*}
\item (normal cone)
\begin{equation*}
(\underline{\bm{x}}-\bm{x})^Tu_x \leq 0, \, \forall\,u_x \in \mathcal{N}_\Omega(\bm{x}), \quad
(\underline{\bm{\lambda}} -\bm{\lambda})^Tu_\lambda \leq 0, \, \forall\,u_\lambda \in \mathcal{N}_{\mathbb{R}^{MN}_+}(\bm{\lambda}),
\end{equation*}
\item (convexity-concavity)
\begin{align*}
&(\underline{\bm{x}}-\bm{x})^Tv_x \leq \tilde{\mathcal{L}}(\underline{\bm{x}},\bm{\lambda}) - \tilde{\mathcal{L}}(\bm{x},\bm{\lambda}), \quad \forall\, v_x\in \partial_{\bm{x}}\tilde{\mathcal{L}}(\bm{x},\bm{\lambda}),\\
&(\underline{\bm{\lambda}}-\bm{\lambda})^Tv_\lambda  \leq \tilde{\mathcal{L}}(\bm{x},\bm{\lambda}) - \tilde{\mathcal{L}}(\bm{x},\underline{\bm{\lambda}}), \quad \forall\, v_\lambda\in \partial_{\bm{\lambda}}(-\tilde{\mathcal{L}})(\bm{x},\bm{\lambda}).
\end{align*}
\end{itemize}
Therefore,
\begin{subequations}\label{eq:condition1}
\begin{align}
(\bm{x}-\underline{\bm{x}})^Tw_x &\leq \tilde{\mathcal{L}}(\underline{\bm{x}},\bm{\lambda}) - \tilde{\mathcal{L}}(\bm{x},\bm{\lambda}), \quad \forall\, w_x\in \Pi_{\Omega}(\bm{x}, - \partial_{\bm{x}}\tilde{\mathcal{L}}(\bm{x},\bm{\lambda})), \\
(\bm{\lambda}-\underline{\bm{\lambda}})^Tw_\lambda & \leq \tilde{\mathcal{L}}(\bm{x},\bm{\lambda}) - \tilde{\mathcal{L}}(\bm{x},\underline{\bm{\lambda}}), \quad \forall\, w_\lambda\in \Pi_{\mathbb{R}^{MN}_+}(\bm{\lambda}, -\partial_{\bm{\lambda}}(-\tilde{\mathcal{L}})(\bm{x},\bm{\lambda}))\text{.}
\end{align}
\end{subequations}
Let $(\bm{x}^*,\bm{\lambda}^*)$ be an equilibrium point of \eqref{eq:distributedAlgorithmCompact}, which satisfies \eqref{eq:equilibrium}. From Lemma \ref{lem:collection} and Theorem \ref{thm:equilibrium}, $\bm{x}^*$ coincides with the unique solution of the primal problem and $(\bm{x}^*,\bm{\lambda}^*)$ is a saddle point of $\tilde{\mathcal{L}}(\bm{x},\bm{\lambda})$. Define
\begin{equation*}
W(\bm{x},\bm{\lambda}) \triangleq \tilde{\mathcal{L}}(\bm{x},\bm{\lambda}^*) - \tilde{\mathcal{L}}(\bm{x}^*,\bm{\lambda}),\, \forall\, \bm{x} \in \Omega, \bm{\lambda} \geq \bm{0}.
\end{equation*}
Obviously, $W(\cdot)$ is a locally Lipschitz continuous function. Moreover, since $(\bm{x}^*, \bm{\lambda}^*)$ is a saddle point of $\tilde{\mathcal{L}}$, $W(\bm{x},\bm{\lambda}) \geq 0$ and $W(\bm{x},\bm{\lambda}) = 0$ if and only if $(\bm{x},\bm{\lambda}) = (\bm{x}^*,\tilde{\bm{\lambda}}^*)$ for some saddle point $(\bm{x}^*,\tilde{\bm{\lambda}}^*)$ of $\tilde{\mathcal{L}}$.

Consider a Lyapunov function
\begin{equation*}
V(\bm{x},\bm{\lambda}) = \frac{1}{2}\|\bm{x} - \bm{x}^*\|^2+\frac{1}{2}\|\bm{\lambda} - \bm{\lambda}^*\|^2.
\end{equation*}
Let $(\bm{x}(t),\bm{\lambda}(t))$ be any solution to \eqref{eq:distributedAlgorithmCompact}. Since $\dot{\bm{x}} \in \mathcal{T}_\Omega(\bm{x})$ and $\dot{\bm{\lambda}} \in \mathcal{T}_{\mathbb{R}^{MN}_+}(\bm{\lambda})$, we have $\bm{x}(t) \in \Omega, \, \bm{\lambda}(t)\geq \bm{0},\, \forall\,t\geq 0$. Moreover, it follows from \eqref{eq:condition1} that
\begin{equation}\label{eq:dV}
\frac{d}{dt}V(\bm{x}(t),\bm{\lambda}(t)) = (\bm{x}(t)-\bm{x}^*)^T\dot{\bm{x}}(t) + (\bm{\lambda}(t)-\bm{\lambda}^*)^T\dot{\bm{\lambda}}(t) \leq  - W(\bm{x}(t),\bm{\lambda}(t)) \leq 0,
\end{equation}
for almost all $t\geq 0$. Therefore, \eqref{eq:distributedAlgorithmCompact} is stable.

Furthermore, since $W(\cdot)$ is locally Lipschitz continuous and $\bm{x}(t),\bm{\lambda}(t)$ are absolutely continuous, $W(t)$ (shorthand for $W(\bm{x}(t),\bm{\lambda}(t))$) is uniformly continuous in $t$.  We claim that $W(t)$ is Riemann integrable over the infinite interval $[0,+\infty)$.   In fact, the Riemann integral of the continuous function $W$ over any finite interval $[0,t)$ equals to the corresponding Lebesgue integral. Moreover, $\int_{0}^tW(\tau)d\tau$ is monotonically increasing since $W$ is nonnegative, and it follows from \eqref{eq:dV} that the Lebesgue integral of $W$ over the infinite interval $[0,+\infty)$, if exists, must be bounded.

As a result, $\int_0^{+\infty}W(\tau)d\tau$ exists and is finite. Then, by the Barbalat's lemma, $(\bm{x}(t),\bm{\lambda}(t))$ converges to the zeros set of $W$, which is exactly the set of saddle points of $\tilde{\mathcal{L}}$.

Let $(\bm{x}^*,\tilde{\bm{\lambda}}^*)$ be a cluster point of $(\bm{x}(t),\bm{\lambda}(t))$ as $t\to +\infty$. Then $(\bm{x}^*,\tilde{\bm{\lambda}}^*)$ is a saddle point of $\tilde{\mathcal{L}}$. Define
\begin{equation*}
\tilde V(\bm{x},\bm{\lambda}) = \frac{1}{2}\|\bm{x} - \bm{x}^*\|^2+\frac{1}{2}\|\bm{\lambda} - \tilde{\bm{\lambda}}^*\|^2.
\end{equation*}
It follows from similar arguments that $\dot {\tilde V }(t)\leq -\tilde{\mathcal{L}}(\bm{x}(t),\tilde{\bm{\lambda}}^*)+ \tilde{\mathcal{L}}(\bm{x}^*,\bm{\lambda}(t))\leq 0$ for almost all $t>0$. Hence, $(\bm{x}^*,\tilde {\bm{\lambda}}^*)$ is Lyapunov stable. It follows from Lemma \ref{lem:semistability} that \eqref{eq:convergence} holds.
\end{IEEEproof}

Finally, we discuss the convergence rate.  Define
\begin{equation}\label{eq:hatxlambda}
\hat{\bm{x}}(t) \triangleq \frac{1}{t}\int_0^t\bm{x}(\tau)d\tau,\quad \hat{\bm{\lambda}}(t) \triangleq \frac{1}{t}\int_0^t\bm{\lambda}(\tau)d\tau\text{,}
\end{equation}
where trajectories $\bm{x}(\cdot),\bm{\lambda}(\cdot)$ are in Theorem \ref{thm:convergence}.

\begin{theorem}\label{thm:rate}%
Under Assumption \ref{assum:1}, there exists a constant $\theta_0>0$ such that
\begin{equation}\label{eq:rate}
\|\tilde{\mathcal{L}}(\hat{\bm{x}}(t),\hat{\bm{\lambda}}(t)) - \tilde{\mathcal{L}}(\bm{x}^*,\tilde{\bm{\lambda}}^*)\| \leq \frac{\theta_0}{t}, \quad \forall\, t>0\text{.}
\end{equation}
\end{theorem}

\begin{IEEEproof}%
Since $\Omega$ is convex and $\bm{x}(\cdot) \in \Omega, \bm{\lambda}(\cdot)\in \mathbb{R}_+^{MN}$, $\hat{\bm{x}}(t)\in \Omega, \hat{\bm{\lambda}}(t)\in \mathbb{R}_+^{MN}, \, \forall\, t>0$. It follows from the Jensen's inequality for the convex-concave $\tilde{\mathcal{L}}$ that, for any $\underline{\bm{x}}\in \Omega, \underline{\bm{\lambda}} \in \mathbb{R}_+^{MN}$,
\begin{equation}\label{eq:Jensen}
\frac{1}{t}\int_0^t \tilde{\mathcal{L}}(\bm{x}(\tau),\underline{\bm{\lambda}})d\tau \geq \tilde{\mathcal{L}}(\hat{\bm{x}}(t),\underline{\bm{\lambda}}), \quad
\frac{1}{t}\int_0^t \tilde{\mathcal{L}}(\underline{\bm{x}},\bm{\lambda}(\tau))d\tau \leq  \tilde{\mathcal{L}}(\underline{\bm{x}},\hat{\bm{\lambda}}(t))\text{.}
\end{equation}
Moreover, it follows from \eqref{eq:condition1} that, for almost all $\tau>0$,
\begin{equation}\label{eq:derivativecondition}
\begin{aligned}
\frac{d}{d\tau}(\frac{1}{2}\|\bm{x}(\tau) - \underline{\bm{x}}\|^2) \leq \tilde{\mathcal{L}}(\underline{\bm{x}},\bm{\lambda}(\tau)) - \tilde{\mathcal{L}}(\bm{x}(\tau),\bm{\lambda}(\tau)),\\
\frac{d}{d\tau}(\frac{1}{2}\|\bm{\lambda}(\tau) - \underline{\bm{\lambda}}\|^2)\leq \tilde{\mathcal{L}}(\bm{x}(\tau),\bm{\lambda}(\tau)) - \tilde{\mathcal{L}}(\bm{x}(\tau),\underline{\bm{\lambda}}).
\end{aligned}
\end{equation}
For any fixed $t>0$, the time average of \eqref{eq:derivativecondition} over integral interval $[0,t]$ with relaxation \eqref{eq:Jensen} indicates
\begin{equation}\label{eq:ratecondition1}
\begin{aligned}
-\frac{\|\bm{x}(0) - \underline{\bm{x}}\|^2}{2t}  \leq \tilde{\mathcal{L}}(\underline{\bm{x}},\hat{\bm{\lambda}}(t)) -\frac{1}{t} \int_0^t \tilde{\mathcal{L}}(\bm{x}(\tau),\bm{\lambda}(\tau))d\tau, \\
-\frac{\|\bm{\lambda}(0) - \underline{\bm{\lambda}}\|^2}{2t} \leq \frac{1}{t} \int_0^t\tilde{\mathcal{L}}(\bm{x}(\tau),\bm{\lambda}(\tau))d\tau - \tilde{\mathcal{L}}(\hat{\bm{x}}(t),\underline{\bm{\lambda}})\text{,}
\end{aligned}
\end{equation}
Replacing $(\underline{\bm{x}}, \underline{\bm{\lambda}})$ by $(\hat{\bm{x}}(t), \hat{\bm{\lambda}}(t))$ in \eqref{eq:ratecondition1} yields
\begin{equation*}
\big\|\frac{1}{t} \int_0^t \tilde{\mathcal{L}}(\bm{x}(\tau),\bm{\lambda}(\tau))d\tau - \tilde{\mathcal{L}}(\hat{\bm{x}}(t),\hat{\bm{\lambda}}(t))\big\|\leq \frac{\theta_1}{2t},
\end{equation*}
where $\theta_1 \triangleq \max_{t>0}\{\|\bm{x}(0) - \hat{\bm{x}}(t)\|^2, \|\bm{\lambda}(0) - \hat{\bm{\lambda}}(t)\|^2\}$. Note that $(\hat{\bm{x}}(t), \hat{\bm{\lambda}}(t))$ are uniformly bounded for $t\in [0,+\infty)$ due to \eqref{eq:convergence}. Similarly, since $\tilde{\mathcal{L}}(\bm{x}^*,\hat{\bm{\lambda}}(t)) \leq \tilde{\mathcal{L}}(\bm{x}^*,\tilde{\bm{\lambda}}^*) \leq \tilde{\mathcal{L}}(\hat{\bm{x}}(t),\tilde{\bm{\lambda}}^*)$, we have from replacing $(\underline{\bm{x}}, \underline{\bm{\lambda}})$ by $(\bm{x}^*,\tilde{\bm{\lambda}}^*)$ in \eqref{eq:ratecondition1} that
\begin{equation*}
\big\|\frac{1}{t} \int_0^t\tilde{\mathcal{L}}(\bm{x}(\tau),\bm{\lambda}(\tau))d\tau - \tilde{\mathcal{L}}(\bm{x}^*,\tilde{\bm{\lambda}}^*) \big\|\leq \frac{\theta_2}{2t},
\end{equation*}
where $\theta_2 \triangleq \max\{\|\bm{x}(0) - \bm{x}^*\|^2, \|\bm{\lambda}(0) - \tilde{\bm{\lambda}}^*\|^2\}$. Thus, \eqref{eq:rate} holds with $\theta_0 \triangleq (\theta_1 + \theta_2)/2$.
\end{IEEEproof}

Theorems \ref{thm:equilibrium}--\ref{thm:rate} provide a complete procedure to prove that algorithm \eqref{eq:distributedAlgorithmCompact} solves problem \eqref{eq:optimizationProblem}. In particular, Theorem \ref{thm:rate} indicates that the value of the Lagrangian function with respect to time average trajectories converges to the value at saddle points with the convergence rate $O(\frac{1}{t})$.

\section{Numerical Examples}
In this section, we first take a simple example for illustration and then consider a more practical example for the performance of our algorithm.
\begin{example}
Consider 4 agents for the optimization problem \eqref{eq:optimizationProblem} with nonsmooth cost and constraint functions:
\begin{equation*}
f_i(x_i) = (x_{i,1} + a_{i,1}x_{i,2})^2 + x_{i,1} + a_{i,2}x_{i,2} + \sqrt{x_{i,1}^2+x_{i,2}^2},\; i = 1,2,3,4
\end{equation*}
and
\begin{equation*}
g_{i,1}(x_i) = \sqrt{x_{i,1}^2 + x_{i,2}^2} - d_{i,1},\quad g_{i,2}(x_i) = -x_{i,1} - x_{i,2} + d_{i,2},
\end{equation*}
where $x_i = (x_{i,1}, x_{i,2}) \in \mathbb{R}^2$, $d_i = (d_{i,1}, d_{i,2}) \in \mathbb{R}^2$ and $a_i = (a_{i,1},a_{i,2})\in \mathbb{R}^2$ for $i = 1,2,3,4$. The local constraint sets of the four agents are
\begin{equation*}
\begin{aligned}
\Omega_1 & = \{x_1 \in \mathbb{R}^2\,|\,(x_{1,1}-2)^2 + (x_{1,2}-3)^2 \leq 25\},\\
\Omega_2 & = \{x_2 \in \mathbb{R}^2\,|\,x_{2,1} \geq 0, x_{2,1} \geq 0, x_{2,1} + 2x_{2,2}\leq 4\},\\
\Omega_3 & = \{x_3 \in \mathbb{R}^2\,|\,4\leq x_{3,1} \leq 6, 2\leq x_{3,2} \leq 5\},\\
\Omega_4 & = \{x_4 \in \mathbb{R}^2\,|\,0 \leq x_{4,1} \leq 15, 0\leq x_{4,1} \leq 20\}.
\end{aligned}
\end{equation*}
The communication graph is shown in Fig. \ref{fig:topology} and algorithm parameters are listed in Table \ref{table:parameters}. Both centralized primal-dual algorithm and our distributed algorithm are utilized to solve this problem and the results are shown in Figs. \ref{fig:Fa}--\ref{fig:Fc}. The trajectories of primal variables are both within their local constraint sets as shown in Figs. \ref{fig:Fa} and \ref{fig:Fb}, while the Lyapunov functions of the algorithms decrease monotonically as shown in Fig. \ref{fig:Fc}.

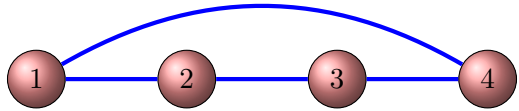
\begin{figure}
\centering
\begin{tikzpicture}[scale = 2]

\tikzstyle{every node}=[shape=circle,draw,minimum size = 20 pt,ball color = red!40]
\path (180 : 1.5cm)   node (1) {$\,1\,$};
\path (180 : 0.5cm)  node (2) {$\,2\,$};
\path (0 : 0.5cm)  node (3) {$\,3\,$};
\path (0 : 1.5cm)  node (4) {$\,4\,$};

\draw [ultra thick] [blue]  (1) -- (2) -- (3) -- (4);
\draw[ultra thick]  [blue] (1) to [out=30,in=150] (4);

\end{tikzpicture}
\caption{The communication graph of the four agents.
}\label{fig:topology}
\end{figure}

\begin{table}
\begin{center}
\caption{Parameters setting}\label{table:parameters}
\begin{tabular}{cccc}
  \hline
   & $d_i$ & $a_i$ & $x_i(0)$ \\ \hline
  $i=1$ & $(6,2)$ & $(8,2)$    & $(2,6)$ \\
  $i=2$ & $(6,3)$ & $(4,7)$    & $(1,1)$ \\
  $i=3$ & $(6,4)$ & $(0.13,8)$ & $(5,4)$ \\
  $i=4$ & $(6,5)$ & $(4,20)$   & $(10,5)$ \\
  \hline
\end{tabular}
\end{center}
\end{table}

\begin{figure}
\begin{center}
\includegraphics[width=15cm]{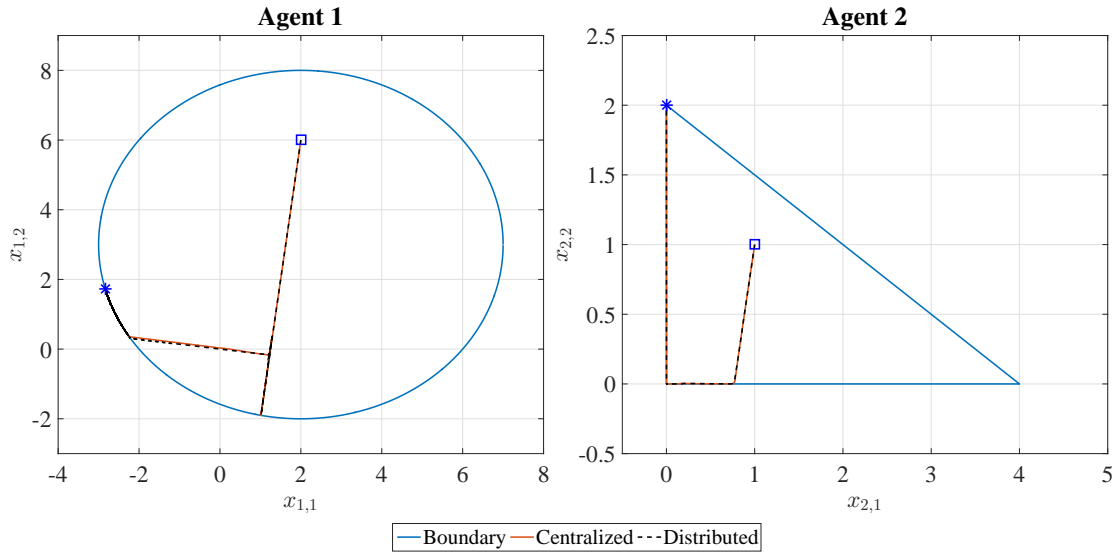}    % The printed column width is 8.4 cm.
\caption{The trajectories of agent 1 and agent 2} \label{fig:Fa}
\end{center}
\end{figure}

\begin{figure}
\begin{center}
\includegraphics[width=15cm]{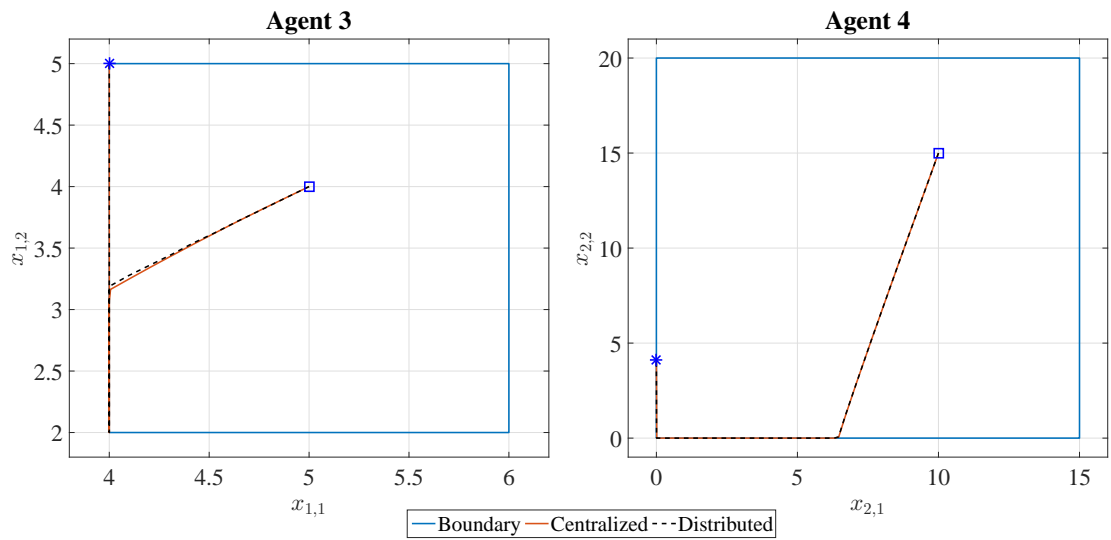}    % The printed column width is 8.4 cm.
\caption{The trajectories of agent 3 and agent 4} \label{fig:Fb}
\end{center}
\end{figure}
\begin{figure}
\begin{center}
\includegraphics[width=15cm]{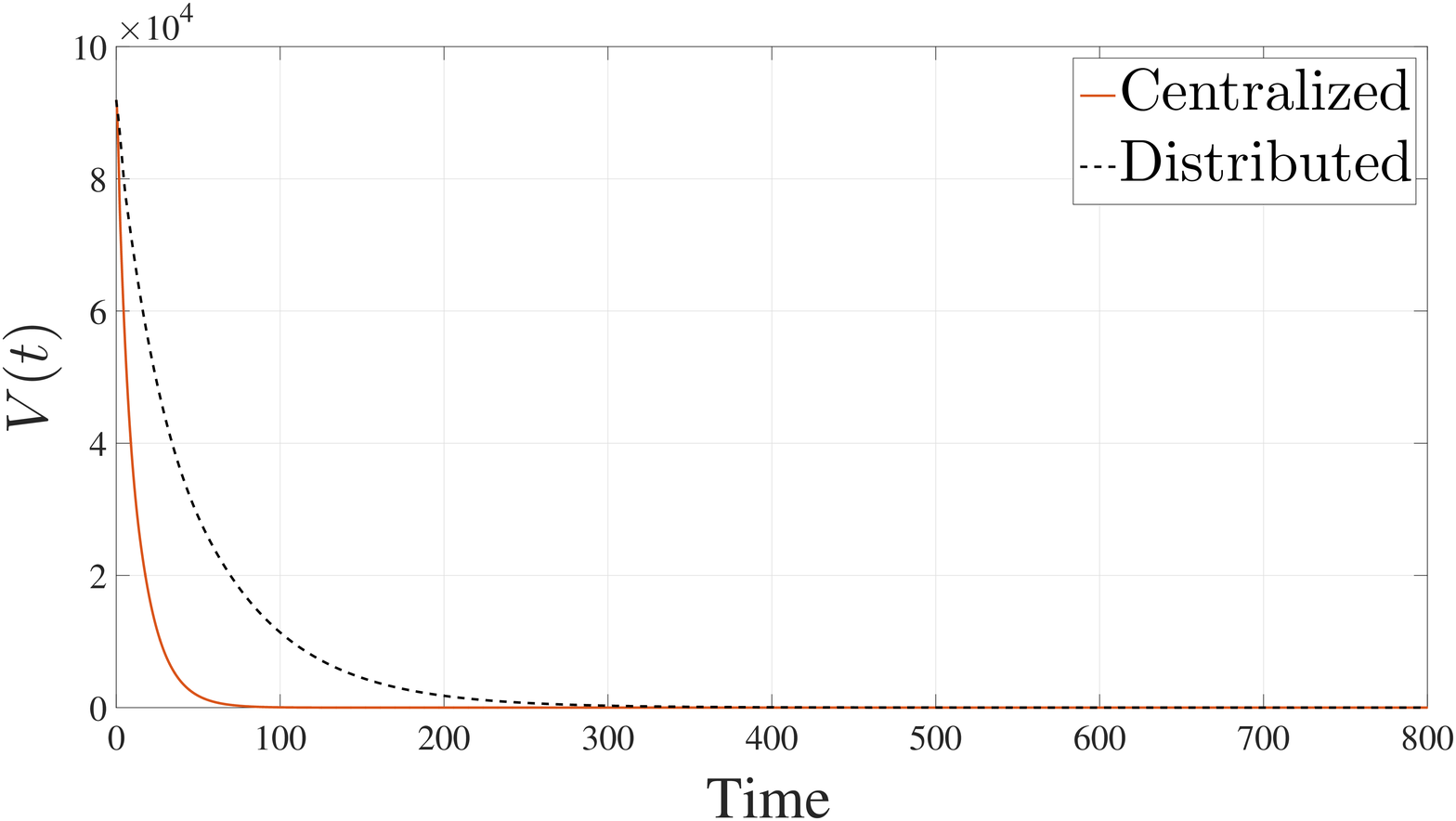}    % The printed column width is 8.4 cm.
\caption{The Lyapunov functions} \label{fig:Fc}
\end{center}
\end{figure}
\end{example}

\begin{example}
Consider problem \eqref{eq:optimizationProblem} with each local constraint as $x_i \in [0,1]$, where each cost function is
\begin{equation*}
f_i(x_i) = a_ix_i^2 + \ln(1+ b_ix_i) + c_i|x_i-d_i| + e_i x_i,
\end{equation*}
and the coupled inequality constraints are $\bm{g}(\bm{x}) = P\bm{x} - q\leq \bm{0}$. We randomly generate coefficients $a_i, b_i, c_i, d_i, e_i \in [0,1]$, matrices $\bm{0}_{M\times N} \leq P \leq \bm{1}_M\bm{1}_N^T$ and $ q \geq \bm{0}_M$ such that strictly feasible point exists. We choose the network size as $N = 10, 20, 50$ and the number of coupled constraints as $M = 5$. For each problem setting, we randomly generate 100 communication graphs. Over each graph, we conduct the numerical experiment and take the relative error $e(t) = \frac{\max_{i =1,...,N}|x_i(t) - x_i^*|}{\max_{i =1,...,N}|x_i^*|}$ for $t = 20, 60, 100$. The average results are shown in Table \ref{tab:ourAlgorithm}, which indicates the effectiveness of our distributed algorithm.
\begin{table}
\begin{center}
\caption{relative error vs. network size}
\label{tab:ourAlgorithm}
\begin{tabular}{cccc}
  \hline
   & $t = 20$ & $t = 60$ & $t = 100$ \\ \hline
  $N=10$ & $0.1982$ & $0.0711$ & $0.0143$ \\
  $N=20$ & $0.5530$ & $0.0290$ & $0.0042$ \\
  $N=50$ & $0.1391$ & $0.0170$ & $0.0105$ \\
  \hline
\end{tabular}
\end{center}
\end{table}
\end{example}
\section{Conclusion}
In this note, a distributed nonsmooth convex optimization problem with coupled inequality constraints has been studied. Based on a modified Lagrangian function constructed via local multipliers and nonsmooth penalty technique, a distributed continuous-time algorithm has been proposed. Also, the convergence of the nonsmooth dynamics has been proved and the convergence rate has been analyzed. Additionally, the effectiveness of the algorithm has been illustrated by two numerical examples.

% or


%
%\appendices
%% you can choose not to have a title for an appendix
%% if you want by leaving the argument blank
%\section{}
%% use section* for acknowledgment
%\section*{Acknowledgment}

\ifCLASSOPTIONcaptionsoff
  \newpage
\fi



\bibliographystyle{IEEEtran}
%\bibliography{refference0,refference1}
\bibliography{E:/HongLab/bib/refference0,E:/HongLab/bib/refference1}
\end{document}